\documentclass[12pt]{amsart} 
\usepackage{amscd}
\usepackage{amssymb}
\usepackage{a4wide}
\usepackage{amstext}
\usepackage{amsthm}
\usepackage{xcolor}
\usepackage{cite}
\usepackage[T1,T2A]{fontenc}
\usepackage[utf8]{inputenc}
\usepackage{url}
\usepackage[T2A]{fontenc}      
\usepackage[russian]{babel}    
\usepackage{amsfonts}
\usepackage{amssymb, amsthm}
\usepackage{amsmath}
\usepackage{mathtools}
\usepackage{needspace}
\usepackage[pdftex]{graphicx}
\usepackage{hyperref}
\usepackage{datetime}
\usepackage{epigraph}
\usepackage{verbatim}
\usepackage{mathtools}
\usepackage{xcolor}
\numberwithin{equation}{section}

\newcommand{\N}{\mathbb{N}}
\newcommand{\R}{\mathbb{R}}

\newcommand{\eps}{\varepsilon}

\newcommand{\CC}{\mathbb {C}}
\newcommand{\Hg}{\mathcal{H}}

\renewcommand{\phi}{\varphi}

\newcommand{\ZZ}{\mathbb{Z}}

\newcommand{\beq}{\begin{equation}}
\newcommand{\eeq}{\end{equation}}
\newcommand{\beqq}{\begin{equation*}}
\newcommand{\eeqq}{\end{equation*}}

\newcommand{\cM}{\mathcal{M}}
\newcommand{\cD}{\mathcal{D}}
\newcommand{\cX}{\mathcal{X}}
\newcommand{\cA}{\mathcal{A}}
\newcommand{\cP}{\mathcal{P}}
\newcommand{\cB}{\mathcal{B}}
\newcommand{\vx}{\vec{x}}
\newcommand{\vy}{\vec{y}}
\newcommand{\ve}{\vec{e}}
\newcommand{\vs}{\vec{s}}
\newcommand{\te}{\tilde{\mathcal E}}
\newcommand{\cE}{{\mathcal E}}

\newtheorem{Thm}{Theorem}[section]
\newtheorem{theorem}[Thm]{Theorem}

\newtheorem{lemma}[Thm]{Lemma}
\newtheorem{proposition}[Thm]{Proposition}
\newtheorem{corollary}[Thm]{Corollary}
\newtheorem{remark}[Thm]{Remark}

\begin{document}
\sloppy
\title[Summation method in optimal control problem with delay]
{Summation method in optimal control problem with delay }
\dedicatory{To N.K.Nikolski on occasion of his 80th birthday}
\author{Pavel Barkhayev}
\address{Pavel Barkhayev, \newline
Department of Mathematical Sciences Norwegian University of Science and Technology,
7491, Trondheim, Norway and
\newline Institute for Low Temperature Physics and Engineering, 47 Nauki ave, 61103, Kharkiv, Ukraine,
\newline {\tt pbarhaev@gmail.com}}

\author{Yurii Lyubarskii}
\address{Yurii Lyubarskii,
\newline St.~Petersburg State University, St. Petersburg, Russia, and
\newline Department of Mathematical Sciences, Norwegian University of Science and Technology, NO-7491 Trondheim, Norway,
\newline {\tt yuralyu@gmail.com} }

\thanks{The work of P.Barkhayev was supported by Norwegian Research Council project No. 275113, "COMAN" 
\newline The work of Yu.Lyubarskii was supported by a grant of the Goverment of the Russian Federation for the state support of scientific research, carried under the supervision of leading scientists, agreement 075-15-2021-602}  \maketitle
 \section{Introduction}
 
  We consider the  control problem for differential equation with  distributed delay
 \beq
 \label{eq:01}
 \dot{x}(t) = x(t-1) + \int_{-1}^0 x(t+\tau) \phi(\tau)d\tau + u(t), \ t > 0
 \eeq
 with the initial conditions 
 \beq
 \label{eq:02} 
 x(t)=x_0(t), \ t\in (-1,0);  \ x(0)= \xi,  \quad  \xi\in \CC,  \ x_0\in L^2(-1,0).
 \eeq
 The distribution of delay is defined by the function $\phi\in L^2(-1,0)$ and the control function $u$ is to be chosen to steer the initial condition $(\xi, x_0(t))$ to the desired state. 
 
 Let 
 $\cM=\CC\times L^2(-1,0)$ be the space of all initial data
 $\cM=\{\vx=\bigl(\begin{smallmatrix}\xi\\ x(t)\end{smallmatrix}\bigr ),  \xi\in \CC,  \ x_0\in L^2(-1,0)\}$.

 The equation (\ref{eq:01}) always has a solution $ x(t)=x(t;\vx,u)$ and 
  we deal with  {\em null controlability} problem: given $T>1$  and $\vx\in \cM$, find the control 
  function $u$ which steers the initial data to zero, i.e.    $ x(t;\vx,u)=0$, $t\in (T-1,T)$.

   It is known (see e.g. \cite{Banks}) that such control  function always exists, and can be given by a bounded operator.
  Actually there exist many such functions.  We consider  the problem of finding the optimal control. In other  
  we look for the problem 
   \beq
   \label{eq:02a}
   \|u\|_{L^2(0,T)} \to \min, \quad u \in U(x_0, \xi;T),
   \eeq
here  $U(x_0,\xi,T)$ stays for all control functions which steers the initial data to zero at the moment $T$.
  
  Such kind of problems were investigated in various settings by many authors. 
  We refer the reader to  (\cite{Halanay}, \cite{Delfour},\cite{Lee}, \cite{Boccia} )  to mention a few. We also refer the reader to \cite{Delfour}, \cite{Boccia} for comprehensive survey 
   and description of the current stay of art. 
   
   We combine the semigroup techniques from \cite{Burns} with the projection techniques, see e.g. \cite{Banksprojections} and 
   \cite{Halebook}. These techniques  lead one to the study of the  eigenfunction expansions  for the corresponding non-selfadjoint operator.
 Since the corresponding systems of eigenfunctions do not form a Riesz basis, the expansions diverge, generally speaking. 
   In \cite{Banksprojections} the authors prove the point-wise convergence on the inner intervals of the segment $[T-1,T]$, yet the problem
   of construction the optimal control on the whole interval remains open. In this article we consider a simplest model of time delay equation and suggest the reconstruction  procedure which converges in the whole $L^2(0,T)$. We believe the same procedure can be applied in much more general settings.
   
 We adjust for this case  the techniques in  \cite{BL}, this  techniques  cannot be applied directly  because   the generating function does not satisfy the Muckenhoupt condition. Also it is worth to mention that when constructing the optimal control for a single eigenvector we arrived
 to already familiar problem
   which appeared in study of the hereditary completeness of sets of exponential functions: given a set of exponents 
  $\Lambda\in \CC$ such that, the corresponding exponential functions $\{\exp(i\lambda t)\}$ are incomplete in $L^2(0,T)$ we need to construct a function which is in the span of  $\{\exp(i\lambda t)\}$ and is orthogonal to all but one 
  of these exponential functions.  This problem has been considered in \cite{BBB} in a very general setting, however the explicit construction
  of such functions is a complicated problem. In the present article we restrict ourselves to the simplest model case when $\phi=0$.

  \section{The homogeneous problem} 
   The corresponding homogeneous problem has the form
   \beqq
 \label{eq:03}
 \dot{x}(t) = x(t-1) + \int_{-1}^0 x(t+\tau) \phi(\tau)d\tau , \ t > 0
 \eeqq
 with the initial conditions 
 \beqq
 \label{eq:04} 
 \bigl(\begin{smallmatrix}x(0)\\ x(\cdot)\end{smallmatrix}\bigr )= \vx_0\in \cM
 \eeqq
 This problem admits the unique solution $x(t)=x(t,\vx_0)$ and also  generates the natural semigroup
 \beqq
 S_t=e^{At}: \cM \to \cM,  \ S_t \vx_0 =  \bigl(\begin{smallmatrix}x(t)\\ x(t+\cdot)\end{smallmatrix}\bigr ).
 \label{eq:05}   
 \eeqq  
 The infinitesimal operator of the group is defined by the relation   
  \beq
  \label{eq:06} 
   A: \vx_0 \mapsto \begin{pmatrix} x_0(-1)+\int_{-1}^0 x_0(\tau)\phi(\tau)d\tau \\
                                                         \dot{x}_0(\tau)
                                \end{pmatrix}.
  \eeq                                                        
   Its domain 
   $$
   \cD_A=\bigl\{(\begin{smallmatrix}\xi \\ x(t)\end{smallmatrix}\bigr ), \ x\in W^2_1(-1,0), \ \xi=x(0)\bigr\}.
   $$  
   
   Let 
   \beqq
   \label{eq:07}
   D(z)=-iz+e^{-iz}+\int_{-1}^0 e^{i\tau z} \phi(\tau) d\tau
   \eeqq
   be the characteristic function of the homogeneous problem. Its indicator diagram is $[0,i]$. We denote by 
     $\Lambda=\{\lambda_n\}$  the set of all zeros of $D$.
  It  can be easily seen that  
   \beqq
   \label{eq:08}
   \lambda_n =-\pi/2+2\pi n+ i \log |n| + o(1), n \in \ZZ.
   \eeqq  
  Here and in what follows we refer the reader to \cite{Levin} for facts about entire functions. We denote by 
 For simplicity we assume $\Lambda \subset \CC_+$ 
  and that $D$ has no multiple zeroes, this may be not true for a finite number of zeros only.
  
 We also have 
 \beq
 \label{eq:08a}
 |D(z)|\asymp \begin{cases}
     |z|, \ \text{for} \ z=x+iy \in \CC_+,  \ y<\log x, \ \text{dist}(z, \Lambda)> \epsilon \\
     e^y,  \ \text{for} \ z=x+iy \in \CC_+,  \ y>\log x, \ \text{dist}(z, \Lambda)> \epsilon. 
                       \end{cases}
\eeq                       
 
 The propositions  below are straightforward (again we refer to \cite{Levin} for the proof of completeness and minimality). 
  
  \begin{proposition}
  The spectra of $A$ is simple and coincides with $i\Lambda$. The corresponding eigenvectors are of the form
$\ve_\lambda=\bigl (\begin{smallmatrix}1\\ \exp(i\lambda \tau)\end{smallmatrix}\bigr ), \lambda \in \Lambda$.
\end{proposition}

Let $\te_\Lambda=\{\ve_\lambda\}_{\lambda \in \Lambda}$ be the set of all eigenfunctions.

\begin{proposition} The system $\te_\Lambda$ is complete and minimal in $\cM$.
\end{proposition}
 
 By $\cX_\Lambda=\{\vx_\lambda\}_{\lambda\in\Lambda}\subset \cM$, $\vx_\lambda= 
 \bigl (\begin{smallmatrix}\xi_\lambda\\ x_\lambda (\tau)\end{smallmatrix}\bigr )$.
 we denote the system biorthogonal to $\te_\Lambda$.
 
 We will use the relations
 \beqq
 \label{eq:08b}
 e^{A\tau}\ve_\lambda = e^{i\lambda\tau}\ve_\lambda, \  e^{A^*\tau}\vx_\lambda = e^{-i\bar{\lambda}\tau}\vx_\lambda,  \
 \lambda \in \Lambda.
 \eeqq
 
 \begin{proposition} The following relations take place
 \beq
 \label{eq:09}
 \overline{\xi_\lambda}= \frac{-i}{D'(\lambda)},
 \eeq
 \beq
 \label{eq:10}
 \int_{-1}^0 e^{izt}\overline{x_\lambda(t)} dt =
 \frac{1}{D'(\lambda)} \frac {e^{-iz}+\int_{-1}^0 e^{iz\tau}\phi(\tau)d\tau +i\lambda}{z-\lambda}.
 \eeq
 \end{proposition}
 
 Since $D(\lambda)=0$, one can  rewrite  (\ref{eq:10} ) as 
 \beq 
 \label{eq:11}
 \int_{-1}^0 e^{izt}\overline{x_\lambda(t)} dt =
  \frac{1}{D'(\lambda)} \bigl [ \frac{e^{-izt}-e^{-i\lambda t}}{z-\lambda} +
                          \int_{-1}^0 \frac{e^{-izt}-e^{-i\lambda t}}{z-\lambda}\phi(t) dt \bigr ].
\eeq
In what follows we   use both   (\ref{eq:10})     and   (\ref{eq:11}).                     
 
 Even after being normalised, the system $\te_\Lambda$ does not form a Riesz basis in $\cM$. However
 for each $\vx\in \cM$ we can consider it formal Fourier series
 \beq
 \label{eq:12}
 \vx \sim \sum_{\lambda\in \Lambda} \langle \vx, \vx_\lambda \rangle \ve_\lambda. 
 \eeq

 \section{Structure of the optimal solution}
 
 
 Denote the operator $\cB: \CC \to \cM$ as $\cB a =  \bigl (\begin{smallmatrix}a\\ 0\end{smallmatrix}\bigr )$, $a\in \CC$.
 The problem (\ref{eq:01}), (\ref{eq:02}) can now be written as
 \beq
 \label{eq:12a}
 \begin{cases}
 \dot{\vx}(t) = A\vx(t)+ \cB u(t), \ t>0, \\
 \vx(0)=\vx_0.
 \end{cases}
 \eeq
  Its solution has the form
  \beqq
  \label{eq:13}
 \vx(t)=e^{At}\vx_0+\int_0^t e^{A(t-\tau)}\cB u(\tau)d\tau.
 \eeqq
 
Fix now $T>1$. The admissibility condition $u\in U(T;\vx_0)$ reads as
 \beqq
 \label{eq:14}
 e^{AT}\vx_0=-\int_0^T e^{A(T-\tau)}\cB u(\tau)d\tau.
 \eeqq
 
 The set of exponential functions $\cE_\Lambda =\{e^{i\lambda t}\}$ is incomplete in $L^2(0,T)$. Let
 $\cE_\Lambda(T)$ be the closure of its span in $L^2(0,T)$ and   
 $\cE_\Lambda(T)^\perp= L^2(0,T)\ominus  \cE_\Lambda(T)$.

We follow the result from \cite{Banksprojections}, see also \cite{Halebook}, Ch.7.4.
  
 \begin{theorem} Let  $u_\lambda \in L^2(0,T)$ solve the problem (\ref{eq:02a}) for $\vx=\ve_\lambda$
 and  \\ $v_\lambda(t)=-u_\lambda(T-t)$. Then
 \beqq
 \label{eq:15}
 v_\lambda\in \cE_\Lambda(T)
 \eeqq
 and also
 \beqq
 \label{eq:16}
 \int_0^T e^{i\mu t} \overline{v_\lambda (t)} dt = \delta_{\lambda, \mu}, \ \mu \in \Lambda.
 \eeqq
  \end{theorem}

 \begin{lemma}
 Let $v\in L^2(0,T)$ be such that 
 \beq
 \label{eq:17}
 0=\int_0^T e^{At } \cB v(t) dt.
 \eeq
 Then $v\in \cE_\Lambda(T)^\perp$.
  \end{lemma}
 
 {\em Proof} Given $\lambda\in \Lambda$ relation (\ref{eq:17}) yields
$$
   0=\int_0^T\langle e^{At } \cB v(t), \vx_\lambda \rangle dt = \int_0^T\langle \cB v(t), e^{A^*t } \vx_\lambda \rangle  dt =
                     \int_0^T\langle \cB v(t), e^{-i\bar{\lambda}t } \vx_\lambda\rangle  dt.
$$
It remains to mention that $ \cB v(t)=\bigl(\begin{smallmatrix}v(t)\\ 0\end{smallmatrix}\bigr )$,
$\vx_\lambda=\bigl (\begin{smallmatrix}\xi_\lambda\\ x_\lambda(t)\end{smallmatrix}\bigr )$,  therefore
$\langle \cB v(t), e^{-i\bar{\lambda}t } \vx_\lambda\rangle =\overline{\xi_\lambda} v(t) e^{-i\bar{\lambda}t}$
and $\xi_\lambda\neq 0$.
\qed
  \begin{lemma}
 Let $\lambda \in \Lambda$ and $v\in L^2(0,T)$ be such that  
 \beqq
 \label{eq:18}
 e^{AT}\ve_\lambda = \int_0^T e^{At } \cB v(t) dt.
 \eeqq
 Then
 \beqq
 \label{eq:19}
 \int_0^T v(t) e^{-i\bar{\mu}t} dt = \bar{y}_\lambda^{-1} e^{i\lambda t}\delta_{\lambda, \mu},  \ \mu \in \Lambda.
 \eeqq
  \end{lemma}
  
  {\em Proof} is similar to that of the previous lemma. We write
  \beqq
  \label{eq:20}
  e^{i\lambda t}\delta_{\lambda, \mu} = \langle e^{AT}\ve_\lambda, \vx_\mu \rangle = \int_0^T \langle e^{At } \cB v(t), \vx_\mu \rangle dt  
  \eeqq
  and then apply   similar arguments. \qed
 
  \begin{corollary}
  The optimal control \underline{linearly} depends on initial data: if $u_1$ and $u_2$ are solutions of the problem (\ref{eq:02a}) for $\vx=\vx_1, \vx_2$ respectively, then $u_1+u_2$ solves this problem for $\vx=\vx_1+\vx_2$.  
  \end{corollary}

Now in order to solve the general problem we need to reconstruct an arbitrary initial data $\vx\in \cM$ from its expansion (\ref{eq:12}).

\section{Summation method}
\subsection{Formulation of the result}

In this section we describe linear summation method for the series   (\ref{eq:12}), i.e we describe the matrices 
$ (w_n(\lambda))_{n\in \N, \lambda \in \Lambda}$ with the following properties:
 \begin{gather}
 \label{eq:21}  w_n(\lambda) \to 1 \ \text{as} \ n\to \infty,  \ \text{for each} \ \lambda\in \Lambda ;\\
 \label{eq:22}  \# \{n; w_n(\lambda)\neq 0\} < \infty, \ \text{for each} \ n\in \N; \\
 \label{eq:23}  S_n \vx \to \vx \ \text{as} \ n\to\infty, \ \text{here} \ S_n\vx=\sum_{\lambda\in \Lambda}
                                      w_n(\lambda)\langle \vx, \vx_\lambda\rangle \ve_\lambda.   
  \end{gather}

  Given the sequences  $\{l_n\}$, $\{R_n\}$, satisfying the conditions
  \begin{gather}
  \label{eq:24}
  l_n, \ R_n \to \infty  \ \text{as} \ n\to \infty,\\
 \label{eq:25}
  l_n^2/n \to 0, \ n/R_n \to 0   \ \text{as} \ n\to \infty,\\
  \label{eq:25a}
  e^{-\pi l_n/2}R_n \to 0  \ \text{as} \ n\to \infty.\
  \end{gather}
  
  \begin{remark} It suffices, for example, to take $l_n= |n|^{1/4}$, $R_n=n^4$. 
  \end{remark}
  
   Consider the function
  \beqq
  \label{eq:26}
  W_n(z)= e^{-l_n\pi- i l_n\log\frac{z-n}{z+n}}, \ z\in \CC_+.
  \eeqq
  
  \begin{theorem} Let, for $\lambda \in \Lambda$
  \label{th:2}
  \beq
  \label{eq:26a}
  w_n(\lambda)=\begin{cases} 
                         W_n(\lambda), & |\lambda|<R_n, \\
                          0, & \text{otherwise}
                           \end{cases} 
 \eeq
 and the operator $S_n$ be defined by 
 \beq
 \label{eq:26b}
 \ S_n\vx=\sum_{\lambda\in \Lambda}
                                      w_n(\lambda)\langle \vx, \vx_\lambda\rangle \ve_\lambda. 
 \eeq
 
 Then $S_n \vx \to \vx$, as $n\to \infty.$

  \end{theorem}
  
  \subsection{Preliminary estimates}

  \begin{lemma}
  \label{l:02}
  \beq
  \label{eq:26c}
  W_n(z) \rightrightarrows 1, \ {\text{as}} \ n\to\infty \ {\text{on each compact set}} \ K\subset \CC_+
  \eeq 
  \end{lemma}
  
  \begin{lemma}
  \label{l:03}
 \beq
 \label{eq:27}
  |W_n(z)|< e^{-\pi l_n/2}, \ z\in \CC_+, |z|=R_n.
  \eeq
  \end{lemma}
  
 \noindent  These lemmas are straightforward.
  
  On the real line we have
  \beqq
  \label{eq:28}
  W_n(t)=
     \begin{cases}
     e^{i  l_n \log \bigl | \frac{t-n}{t+n} \bigr |}, & |t|<n, \\
       e^{i l_n \log \bigl | \frac{t-n}{t+n} \bigr | d -l_n\pi}, & |t|> n
     \end{cases}
  \eeqq

  \begin{lemma}
  \label{l:04}
  \beqq
  \label{eq:29}
  \sup_n \bigl |\int_{1<|t|<n} W_n(t) \frac{dt}t\bigr | < \infty
 \eeqq 
 \end{lemma}
  
  \begin{lemma}
  \label{l:05}
  \beqq
  \label{eq:30}
  \sup_n \bigl |\int_{|t|>n} W_n(t) \frac{dt}t\bigr | < \infty
 \eeqq
 \end{lemma}
 
 \bigskip
 
 \noindent {\em Proof of Lemma \ref{l:04}}
 
 We have 
 \beqq
 \label{eq:31}
 I_n:= \int_{1<|t|<n} W_n(t) \frac{dt}t = \bigl (\int_{1<|t|<n/l_n} + \int_{n/l_n<|t|<n} \bigr ) = I_{n,1}+I_{n,2}
 \eeqq
  
  The estimate of $I_{n,1}$ is straightforward. We have for $|t|<n/l_n$ 
   $$
 l_n  \log \bigl |\frac{t-n}{t+n}| \prec \frac {l_n}n |t| << 1,
   $$
 Respectively
 $$
  I_{n,1}= \int_{1<|t|<n/l_n} \bigl [1+ O\bigl (\frac {l_n}n |t| \bigr ) \bigr] \frac{dt}{t} = O(1).
  $$
  
 \medskip
 
When estimating $I_{2,n}$ we consider   each interval $\pm [n/l_n, n] $ separately.
Consider the integral
$$
I_{n,2}'=  \int_{n/l_n<t<n}      e^{i  l_n \log \bigl | \frac{t-n}{t+n} \bigr | }\frac{dt}t
$$
 By changing of variables 
 $$
u= \frac{n-t}{n+t}; \ \tau=1-u; \  v= \log(1-\tau),
$$
and extracting the main term on each step, we finally arrive to the integral
$
\int_{1/l_n}^\infty e^{il_n v}\frac{dv}v,
$  
 which is bounded uniformly in $n$. \qed
 
 \medskip
 
 {\em Proof of Lemma \ref{l:05}} goes in the same vien. We  
 write
 $$
 \int_{|t|>n} W_n(t) \frac{dt}t= \bigl ( \int_{|t|>nl_n} + \int_{n<|t|<nl_n} \bigr ) W_n(t) \frac{dt}t
 $$
 and use the similar reasonings. \qed
 
 \subsection{Proof of Theorem \ref{th:2}}
 {\bf a.} Let operator $S_n$ be defined by (\ref{eq:23}). It follows from (\ref{eq:26c}) that 
 $S_n\vx \to \vx$ as $n\to \infty$ on the  dense set of finite linear combinations of $\ve_\lambda$. 
 So it suffices to prove the uniform boundedness of all operators $S_n$.   
Consider instead the adjoint operator
 \beqq
 \label{eq:32}
 S^*\vy=\sum_\lambda \overline {w_n(\lambda)}\langle \vy, \ve_\lambda\rangle \vx_\lambda.
 \eeqq
 Actually it is more convenient to deal with the conjugate objects:
 \beqq
 \label{eq:33}
\overline{ S^*\vy}=\sum_\lambda {w_n(\lambda)}\langle  \ve_\lambda, \vy\rangle \overline{\vx_\lambda} = :\vs_n.
\eeqq 
 Let $\vs_n=\bigl (\begin{smallmatrix}\sigma_n\\ s_n(t)\end{smallmatrix}\bigr )$ and  $\vy=\bigl (\begin{smallmatrix}\eta\\y(t)\end{smallmatrix}\bigr )$.  
  Since $\vx_\lambda=\bigl (\begin{smallmatrix}\xi_\lambda\\ x_\lambda(t)\end{smallmatrix}\bigr )$ we have
  \beq
  \label{eq:34}
 \sigma_n=\sum_\lambda w_n(\lambda) \bar{\xi}_\lambda + 
   \sum_\lambda w_n(\lambda) \bar{\xi}_\lambda \int_{-1}^0e^{i\lambda t} \overline {y(t)}dt =\sigma_n^{(1)}+
   \sigma_n{(2)},
 \eeq  
 \beq
 \label{eq:35}
 s_n(t)=\sum_\lambda w_n(\lambda)\overline{x_\lambda(t)} + \sum_\lambda w_n(\lambda)\int_{-1}^0e^{i\lambda \tau}
 \overline{y(\tau)} d\tau \ \overline{x_\lambda(t)}= s_n^{(1)}(t)+ s_n^{(2)}(t).
 \eeq 
 Each of the four summands in the right of (\ref{eq:34}), (\ref{eq:35}) should be estimates separately. 
 
 Consider the contour $\Gamma_n=L_n\cup C_n$. Here $L_n=[-R_n,R_n]$, $C_n=\{R_n e^{i\theta}; \theta\in [0,\pi]\}$. 
The arcs $C_n$ uniformly in $n$ satisfy the Carleson condition, hence, for any
 function $h$ from the Hardy space $\Hg(\CC_+)$,
 \beqq
  \label{eq:35a}
  \Bigl ( \int_{C_n} |h(\zeta)|^2|d\zeta| \Bigr )^{1/2} \leq {\rm Const} \|h\|.
  \eeqq
 We  also assume that ${\rm dist}(\Lambda, \Gamma_n)>0$ uniformly in $n$. 

 \medskip
 
 {\bf b.} \underline {Estimate of $\sigma_n^{(1)}$}. We use (\ref{eq:09}). By the residues 
 \beq
 \label{eq:36}
 2i\pi \sigma_n^{(1)} = \int_{\Gamma_n} \frac {W_n(\zeta)}{D(\zeta)}d\zeta = \bigl (\int_{C_n}+\int_{L_n}\bigr ) \frac {W_n(\zeta)}{D(\zeta)}d\zeta = I_n + J_n.
\eeq 

That $I_n$ is bounded (and even tends to zero as $n\to 0$) is very clear, this follows from (\ref{eq:25a}), (\ref{eq:27}), and (\ref{eq:08a}). In order to see the boundedness of $J_n$ it suffices to mention that, according Lemma  \ref{l:05} the integral
$$
\int_{L_n} \frac{W_n(t)}{it+1}dt
$$
is bounded uniformly in $n$. Therefore it suffices to establish the boundedness of
 $$
\int_{L_n}W_n(t) \bigl \{ \frac 1{D(t)}-\frac 1 {it+1}\bigr\} dt.
$$
The later is evident since the expression in parenthesis  uniformly $O(t^{-2})$ as $t\to \infty$.

\medskip

{\bf c.} \underline{Estimate of $\sigma_n^{(2)}$} Denote 
\beq
\label{eq:36a}
Y(z)=\int_{-1}^0 e^{iz\tau} \overline{y(\tau)}d\tau.
\eeq
As before we have 
\beqq
\label{eq:37}
 2i\pi \sigma_n^{(2)} = \int_{\Gamma_n}  W_n(\zeta)\frac{Y(\zeta)}{D(\zeta)}d\zeta = \Bigl (\int_{C_n}+\int_{L_n}\Bigr )  W_n(\zeta)\frac{Y(\zeta)}{D(\zeta)}  d\zeta = I_n + J_n.
\eeqq 
The estimate of $J_n$ is now straightforward because $Y\in L^2(\R)$ uniformly in $y$ in the unit ball of $L^2(-1,0)$ and
$|D(t)|\asymp |t|$ as $t\to \pm \infty$.

In order to estimate $I_n$ we observe that the function $\Psi(\zeta)=e^{iz}Y(\zeta)$ belongs to the Hardy space $H^2(\CC_+)$
and has the same $L^2(\R)$ norm as $Y$. Also, it follows from (\ref{eq:08a}) that 
$$
\Bigl |\frac{Y(\zeta)}{D(\zeta)} \Bigr | \prec |\Psi(\zeta)|, \ \zeta \in C_n.
$$
 Now we use that  all $C_n$'s  are uniformly  Carleson curves  in $\CC_+$ and also (\ref{eq:27}). We obtain
 $$
 [J_n| \prec e^{-l_n\pi/2} \int_{C_n} |\Psi(\zeta)||d\zeta | \leq   
         e^{-l_n\pi/2} |C_n|^{1/2} \bigl (\int_{C_n} |\Psi(\zeta)|^2 |d\zeta|\bigr )^{1/2}
 \prec R_n^{1/2} e^{-l_n\pi/2} \|\Psi\|.
 $$
 It remains to refer to (\ref{eq:25a}).

\medskip

{\bf d.} \underline{Estimate of $\|s_n^{(1)}\|$.} We apply the Fourier transform and use representation 
(\ref{eq:11}):
\beqq
\label{eq:38}
({s_n^{(1)}})\hat{}(x)=\sum_\lambda w_n(\lambda)\frac 1 {D'(\lambda)}
 \Bigl (\frac{e^{-izt}-e^{-i\lambda t}}{z-\lambda} +
                          \int_{-1}^0 \frac{e^{-izt}-e^{-i\lambda t}}{z-\lambda}\phi(t) dt \Bigr )
\eeqq
Therefore it suffices to obtain uniform in $t$ estimates of $\|\Phi_{n,t}\|_{L^2(\R)}$
where 
\beqq
\label{eq:39} 
\Phi_{n,t}(x)=\sum_\Lambda w_n(\lambda) \frac 1 {D'(\lambda)}
       \frac{e^{-ixt}-e^{-i\lambda t}}{x-\lambda}
\eeqq
Since $\Lambda \subset \CC_+$ this function belongs to the Hardy space $\Hg^2(\CC_-)$
and 
\beqq
\label{eq:40}
\|\Phi_{n,t}\|_{L^2(\R)}=\sup \Bigl \{\bigl | \int_{-\infty}^\infty \Phi_{n,t}(x)h(x)dx \Big | \ ; \ 
h\in \Hg^2(\CC_+), \ \|h\| = 1 \Bigr \}.
\eeqq
We obtain two terms to be estimated:
\beqq
\label{eq:41}
\sum_\Lambda  \frac {w_n(\lambda)}{D'(\lambda)}e^{i\lambda t}\int_{-\infty}^\infty\frac{h(x)}{x-\lambda}dx=
2i\pi \sum_\Lambda  \frac {w_n(\lambda)}{D'(\lambda)}e^{i\lambda t}h(\lambda)=:a_{1,n},
\eeqq
and
\beq
\label{eq:42}
\sum_\Lambda  \frac {w_n(\lambda)}{D'(\lambda)}\int_{-\infty}^\infty\frac{e^{i x t}h(x)}{x-\lambda}dx=:a_{2,n}.
\eeq

Once again we have 
\beqq
\label{eq:43}
a_{1,n}=\Bigl (\int_{L_n}+\int_{C_n}\Bigr )W_n(\zeta)\frac {e^{i\zeta t}}{D(\zeta)}h(\zeta)d\zeta.
\eeqq
The estimate of the first integral follows blantly from. the Schwartz inequality. In order to estimate the second one we observe that $\bigl |e^{it \zeta}/D(\zeta) \bigr | \prec 1$, $\zeta \in C_n$ uniformly in $t$ and $n$ and then again use the fact that all $C_n$'s are uniformly Carlesonian together with (\ref{eq:25a}):
\beqq
\label{eq:44}
\Bigl | \int_{C_n} W_n(\zeta)\frac {e^{i\zeta t}}{D(\zeta)}h(\zeta)d\zeta \Bigr | \prec 
        \int_{C_n} |W_n(\zeta)| |h(\zeta)| |d\zeta| \prec e^{-\pi l_n/2} R_n^{1/2} \|h\|\prec \|h\|.
        \eeqq
        
 In order to estimate $a_{n,2}$ we represent $h$ as
 $$
 h(x)=\int_0^\infty e^{i\xi x}\omega(\xi)d\xi
 $$
 for some $\omega \in L^2(0, \infty)$.
 
 For $t\in [-1,0]$ we have
\beq
\label{eq:44a}
 h(x)=\int_t^0 e^{i\xi x}\omega(\xi-t)d\xi +\int_0^\infty e^{i\xi x}\omega(\xi-t)d\xi= h_t^-(x)+h_t^+(x).
 \eeq
Then $h_t^\pm \in \Hg(\CC_\pm)$, $\|h_t^+\| \leq \|h\|$ and 
\beqq
\label{eq:45}
a_{2,n}=\sum_\lambda w_n(\lambda)\frac{h_t^+(\lambda)}{D'(\lambda)}.
\eeqq
The rest of the estimate is similar to that of $a_{n,1}$.

\medskip

{\bf e.} \underline{Estimate of $\|s_n^{(2)}\|$.}  Again we apply the Fourier transform but now use representation 
(\ref{eq:10}):
\beqq
\label{eq:46}
(\widehat{s_n^{(2)}})(x)=\sum_\lambda w_n(\lambda)\frac {Y(\lambda)} {D'(\lambda)} \ 
 \frac{ e^{-ix}+\int_{-1}^0e^{ix\tau}\phi(\tau)d\tau+i\lambda}{x-\lambda}, 
\eeqq
here $Y(\lambda)$ is defined in (\ref{eq:36a}).

Once again we use the duality reasons
$$
\|\widehat{s_n^{(2)}}\|= \sup \bigl \{ \bigl | \int_{-\infty}^\infty \widehat{s_n^{(2)}}(x) h(x)dx \bigr | \ , \ h\in\Hg^2(\CC_+), \
\|h\|\leq 1 \bigr \}.
$$
Let the functions $h_t$ be defined as in (\ref{eq:44a}). Consider the function
$$
f(x)=h_{-1}+\int_{-1}^0h_t(x)\phi(t)dt.
$$
We have $f\in \Hg^2(\CC^+)$, $\|f\| \prec 1$ and 
\beqq
\label{eq:47}
 \int_{-\infty}^\infty \widehat{s_n^{(2)}}(x) h(x)dx= \sum_\lambda w_n(\lambda)\frac{1+i\lambda}{D'(\lambda)}Y(\lambda)
 f(\lambda).
 \eeqq
 This expression can be estimated in a similar way as when estimating $\|s_{n,1}\|$.
  \qed
 
 We have proved the following statement.
 
 \begin{theorem}
 \label{th:3}
 Given the control problem (\ref{eq:01}), (\ref{eq:02}), let $u_0(t)$ be the optimal solution defined by 
 (\ref{eq:02a}). Let, further,  $\Lambda$ and $\te_\Lambda$ be the corresponding eigenvalues and eigenvectors
 and also 
 \beqq
 \vx_0 \sim \sum_{\lambda\in \Lambda} \langle \vx_0, \vx_\lambda \rangle \ve_\lambda
 \eeqq
be the formal Fourier expansion of the initial data $\vx_0$.

Then  
\beqq
\label{eq:47a}
u_0(t)=\lim_{n\to \infty} \sum_\lambda w_n(\lambda) \langle \vx_0, \vx_\lambda \rangle u_\lambda(t), 
\eeqq
here $w_n(\lambda)$ are defined in (\ref{eq:26a}) and the functions $u_\lambda \in L^2(0,T)$ are defined by the 
relations
\beqq
\label{eq:47b}
u_\lambda(t)= -\overline{v_\lambda(T-t)}; \  v_\lambda\in  \cE(\Lambda), \ 
\int_0^T e^{i\mu t} \overline{v_\lambda(t)}dt = \delta_{\lambda, \mu}, \ \mu\in\Lambda.
\eeqq
The limit exists in $L^2(0,T)$ sense.
\end{theorem}

  \section{A model case}
 
 {\bf a.}  In this section we consider the simplest model case
 \begin{gather}
 \label{eq:48}
 \dot{x}(t)=x(t-1)+u(t), \ t >0;\\
 \label{eq:49}
 x(t)=x_0(t), \ t\in (-1,0), \ x(0)=\xi_0.
 \end{gather}
 For this case we give an explicit construction of the corresponding functions $v_\lambda$ for $T\in (1,2)$.
 We denote  $\delta=T-1$.
 The characteristic  function has the form
 \beq
 \label{eq:49a}
 D(z)=-iz+e^{-iz}. 
 \eeq
 
 \medskip
 
 {\bf b.} Description of $\cE(\Lambda)^\perp$. Let $g\in L^2(0,T)$, $g\in \cE(\Lambda)^\perp$. Then the entire function
 \beqq
 \label{eq:50}
 G(z)=\int_0^T e^{izt} \overline{g(t)}dt
 \eeqq
 vanishes on $\Lambda$ and the standard reasonings on comparing its growth with that of $D(z)$ yield
 \beq
 \label{eq:51}
 G(z)= e^{iz} D(z) \Omega(z),
 \eeq
 where $\Omega$ runs through the set of all functions of the form
 \beq
 \label{eq:52}
 \Omega(z)=\int_0^\delta e^{izt} \overline{\omega(t)}dt,  \ x\Omega(x)\in L^2(\R)
 \eeq
 The later yields
 \beq
 \label{eq:53}
 \omega\in W^2_1(0,\delta), \   \omega(0)=\omega(\delta)=0,
 \eeq
 and, by (\ref{eq:49a}) - (\ref{eq:51}), we obtain
 \beqq
 \label{eq:54}
 g(t)=\begin{cases}
                       \omega'(t-1),  & t\in (1,1+\delta), \\
                       \omega(t), & t\in (0,\delta), \\
                       0, & \rm{otherwise}
         \end{cases} 
 \eeqq
 This is the description of $\cE(\Lambda)^\perp$. 
 
 \medskip
 
 {\bf c.} Description of the biorthogonal system. Fix $\lambda\in \Lambda$ and let a function 
 $v_\lambda\in L^2(0,T)$ satisfy
 \beqq
 \label{eq:55}
  \int_0^T e^{i\mu t} \overline{v_\lambda(t)}dt = \delta_{\lambda, \mu}, \ \mu\in\Lambda.
 \eeqq
 Denote
 \beqq
 \label{eq:56}
V_\lambda(z)=\int_0^T e^{iz t} \overline{v_\lambda(t)}dt . 
 \eeqq
 As before we have
 \beqq
 \label{eq:57}
  V_\lambda(z)=\frac 1 {e^{i\lambda}D'(\lambda)}\frac{e^{iz}D(z)}{z-\lambda}P_\lambda(z),
 \eeqq
 where the function
 \beqq
 \label{eq:58}
 P_\lambda(z)=\int_0^\delta e^{iz t} \overline{p_\lambda(t)}dt,  \ p_\lambda\in L^2(0,\delta)
 \eeqq
 Satisfies
 \beq
 \label{eq:59}
 P_\lambda(\lambda)=\int_0^\delta e^{i\lambda t} \overline{p_\lambda(t)}dt=1.
 \eeq
 
 We remind that $D(\lambda)=0$. Therefore 
 \beqq
 \label{eq:60} 
 \frac{e^{ix}D(x)}{x-\lambda}=-ie^{ix}-i \frac{e^{i(x-\lambda)}-1}{x-\lambda}=
        -i \left [ e^{ix}+\int_0^1 e^{ixt} e^{- i\lambda t} dt \right ]
 \eeqq       
 and, denoting $a_\lambda=1/(e^{i\lambda} D'(\lambda)$ we obtain
 \beqq
 \label{eq:61}
 v_\lambda(t) = -ia_\lambda \left [\overline{p_\lambda(t-1)}+ \int_{t-1}^te^{-i\lambda(t-u)}\overline{p_\lambda(u)} du 
                                              \right ]
    \eeqq                                          
When  $p_\lambda$ runs through $L^2(0,\delta)$ the function $v_\lambda$ runs through all functions in $L^2(0,T)$
orthogonal to $\cE_\Lambda \setminus \{e^{i\lambda t}\}$. The condition $v_\lambda\in \cE_\lambda^\perp$ leads one to condition on  $p_\lambda$.

\medskip

{\bf d.} We assume that the functions $\omega$, $p_\lambda$ are defined on the whole $\R$ and vanish outside of 
$(0, \delta)$. The condition $v_\lambda \in \cE_\Lambda$ takes the form

\beq
\label{eq:62}
0=\int_0^T [\omega(t)+\omega'(t-1)]  \left [\overline{p_\lambda(t-1)}+ e^{- i\lambda t}
\int_{t-1}^te^{-i\lambda \tau}\overline{p_\lambda(\tau)} d\tau \right ]dt,  
\eeq
 for all $\omega\in W^2_1(0, \delta)$, $\omega(0)=\omega(\delta)=0$.
 
 When opening the parenthesis we obtain four summands $A_i$, $i=1,2,3,4$. Consider each of them separately.
 \beqq
 \label{eq:63}
 A_1= \int_0^T \omega(t)\overline{p_\lambda(t-1)}dt=0,
 \eeqq
 since the supports of $\omega(t)$ and $p_\lambda(t-1)$ are disjoint.
 
 Denote $\eps_\lambda(t)=e^{-i\lambda t}\chi_{(0, \infty)}(t)$. Then 
 \beqq
 \label{eq:64}
 A_2=\int_0^T\omega(t) \int_{t-1}^t    \overline{p_\lambda(\tau)}e^{-i\lambda(t-\tau)}d\tau=
 \int_0^\delta\omega(t)( \overline{p_\lambda} * \eps_\lambda) (t)dt;
 \eeqq
 \beqq
 A_3=\int_0^T\omega'(t-1)\overline{p_\lambda(t-1)}dt = -
                   \int_0^T\omega(t-1)\overline{p_\lambda'(t-1)}dt.
 \label{eq:65}
 \eeqq
 Here we assume that $p_\lambda'$ exists. With this assumption we will arrive to equation with respect to $p_\lambda'$ 
 which has smooth solution. This suffices since we already know that there is only one $p_\lambda$, satisfying all conditions.
 
 When considering $A_4$ we use (\ref{eq:59}). We have
 \begin{gather*}
 \label{eq:66}
 A_4=\int_0^T \omega'(t-1)e^{-i\lambda t}\int_{t-1}^t \overline{p_\lambda(\tau)}e^{i\lambda\tau}d\tau dt=
 \int_1^T \omega'(t-1)e^{-i\lambda t}\int_{t-1}^\delta\overline{p_\lambda(\tau)}e^{i\lambda\tau}d\tau dt= \\
 e^{-i\lambda}\int_0^\delta\omega'(t)e^{-i\lambda t} \left (
           1-\int_0^t\overline{p_\lambda(\tau)}e^{i\lambda\tau}d\tau 
                                                                                  \right )dt=
 e^{-i\lambda} \left [i\lambda      \int_0^\delta\omega(t)e^{-i\lambda t} dt -
            \int_0^\delta  \omega(t)     ( \overline{p_\lambda} * \eps_\lambda) '(t) dt \right ]   .                                                                     
 \end{gather*}
 
 We collecting these relations together and remark that $\omega$ can be any function in $W^2_1(0,\delta)$, satisfying
 (\ref{eq:53}). We then arrive to the  equation  
 \beq
 \label{eq:67}
 -\overline{p_\lambda}'(t)-e^{-i\lambda} ( \overline{p_\lambda} * \eps_\lambda) '(t) +( \overline{p_\lambda} * \eps_\lambda) (t)
 - e^{-i\lambda}\eps_\lambda'(t)=0,  \ 0<t<\delta.
 \eeq
The left-hand side of this equation is well-defined for all $t>0$. Denote it by $\alpha(t)$.  
Let $\cP_\lambda(s)$ and $e^{-\delta s} \cA(s)$ be the Laplace transform of $\overline{p_\lambda}$ and $\alpha$ respectively.
Taking once again into account  that $D(\lambda)=0$, we transform (\ref{eq:67}) into equation for   $\cP_\lambda(s)$ :
 \beqq
 \label{eq:68}
 \cP_\lambda(s)(1-s^2)= (\lambda^2-i\lambda \overline{p_\lambda(0)}) +s\overline{p_\lambda(0)} + e^{-\delta s} \cA(s).
 \eeqq
 
 The last summand in the right-hand side does not influence the values of $ \overline{p_\lambda}$  on $(0,\delta)$.
 Therefore 
 \beqq
 \label{eq:69}
 \overline{p_\lambda}(t)= [\lambda^2-i\lambda \overline{p_\lambda}(0)]\sinh t + \overline{p_\lambda}(0)\cosh t, \ t\in (0,\delta).
 \eeqq

 The value $\overline{p_\lambda}(0)$ should now be chosen to meet (\ref{eq:59}). We omit the corresponding calculation.

\end{document}